\documentclass[11pt]{amsart} 
\usepackage{verbatim, latexsym, amssymb, amsmath,color}
\usepackage{epsfig}
\usepackage{tikz}
\usepackage{hyperref}
\usepackage{float}
\usepackage{graphicx}
\usepackage{sidecap}

\usetikzlibrary{matrix}

\def\R{\mathbb R}
\def\N{\mathbb N}
\def\C{\mathcal C}
\def\H{\mathcal H}
\def\T{\mathcal T}
\def\pprime{\prime\prime}

\theoremstyle{remark}
\theoremstyle{definition}

\title[Width and positive scalar curvature]{Metrics of positive scalar curvature and unbounded widths}
\author{Rafael Montezuma}
\address{Instituto de Matem\'atica Pura e Aplicada (IMPA) \\ Estrada Dona Castorina 110 \\ 22460-320 Rio de Janeiro \\ Brazil}
\email{rafaelmc@impa.br}

\thanks{The author was partly supported by FAPERJ and NSF}

\begin{document}

\begin{abstract}
{ {In this work we construct a sequence of Riemannian metrics on the three-sphere with scalar curvature greater than or equal to $6$ and arbitrarily large widths. Our procedure is based on the connected sum construction of positive scalar curvature metrics due to Gromov and Lawson. We develop analogies between the area of boundaries of special open subsets in our three-manifolds and $2$-colorings of associated full binary trees. Then, via combinatorial arguments and using the relative isoperimetric inequality, we argue that the widths converge to infinity.}}
\end{abstract}
\maketitle
\setcounter{tocdepth}{1}
%\tableofcontents

\section{Introduction}\label{introduction}

Since the proof of the positive mass conjecture in general relativity by Schoen and Yau \cite{S-Y}, and Witten \cite{Witten}, the rigidity phenomena involving the scalar curvature has been fascinating the geometers. These results play an important role in modern differential geometry and there is a vast literature about it, see (\cite{Ambrozio, Bray, BBEN, BBN, B-M, BMN, C-G, Eichmair, M-N, Miao, M-M, Min-Oo, Moraru, SY}). Many of these works concern rigidity phenomena involving the scalar curvature and the area of minimal surfaces of some kind in three-manifolds.

The width of a Riemannian three-manifold $(M^3, g)$ is a very interesting geometrical invariant which is closely related to the production of unstable, closed, embedded minimal surfaces, see \cite{C-D} or \cite{Pitts}. It can be defined in different ways depending on the setting, but has always the intent to be the lowest value $W$ for which it is possible to sweep $M$ out using surfaces of $g$ area at most $W$. Let us be more precise about this object. 

%In this work we deal with the width
Let $g$ be a Riemannian metric on the three-sphere. A \textit{sweepout} of $(S^3, g)$ is a one-parameter family $\{\Sigma_t\}$, $t\in [0,1]$, of smooth $2$-spheres of finite area which are boundaries of open subsets $\Sigma_t = \partial \Omega_t$, vary smoothly and degenerate to points at times zero, $\Omega_0 = \varnothing$, and one, $\Omega_1 = S^3$. 

The simplest way to sweep out the three-sphere is using the level sets of any coordinate function $x_i : S^3\subset \R^4 \rightarrow \R$. 

Let $\Lambda$ be a set of sweepouts of $(S^3, g)$. It is said to be saturated if given a map $\phi \in C^{\infty}([0,1]\times S^3, S^3)$ such that $\phi(t, \cdot)$ are diffeomorphism of $S^3$, all of which isotopic to the identity, and a sweepout $\{\Sigma_t\} \in \Lambda$, we have $\{\phi(t,\Sigma_t)\} \in \Lambda$. 

The \textit{width} of the Riemannian metric $g$ on $S^3$ with respect to the saturated set of sweepouts $\Lambda$ is defined as the following min-max invariant:
\begin{equation*}
W(S^3, g) = \inf_{\{\Sigma_t\} \in \Lambda} \ \ \sup_{t\in [0,1]} Area_g(\Sigma_t),
\end{equation*}
where $Area_g(\Sigma_t)$ denotes the surface area of the slice $\Sigma_t$ with respect to $g$.

%M-N rigidez da width
Marques and Neves \cite{M-N} proved that the width of a metric $g$ of positive Ricci curvature in $S^3$, with scalar curvature $R\geq 6$, satisfies the upper bound $W(S^3, g) \leq 4\pi$ and there exists an embedded minimal sphere $\Sigma$, of index one and surface area $Area_g(\Sigma) = W(S^3, g)$. They proved also that in case of equality $W(S^3, g) = 4\pi$, the metric $g$ has constant sectional curvature one. The main purpose of the present work is to prove that this is no longer true without the assumption on the Ricci curvature. More precisely, we have:

\subsection*{Theorem A}\label{thm.A}
\textit{For any $m>0$ there exists a Riemannian metric $g$ on $S^3$, with scalar curvature $R\geq 6$ and width $W(S^3, g) \geq m$.
}

\subsubsection*{Remark} It is interesting to stress that the Riemann curvature tensors of the examples that we construct are uniformly bounded. 

\begin{figure}[H]
  \centering
    \includegraphics[width=.75\textwidth]{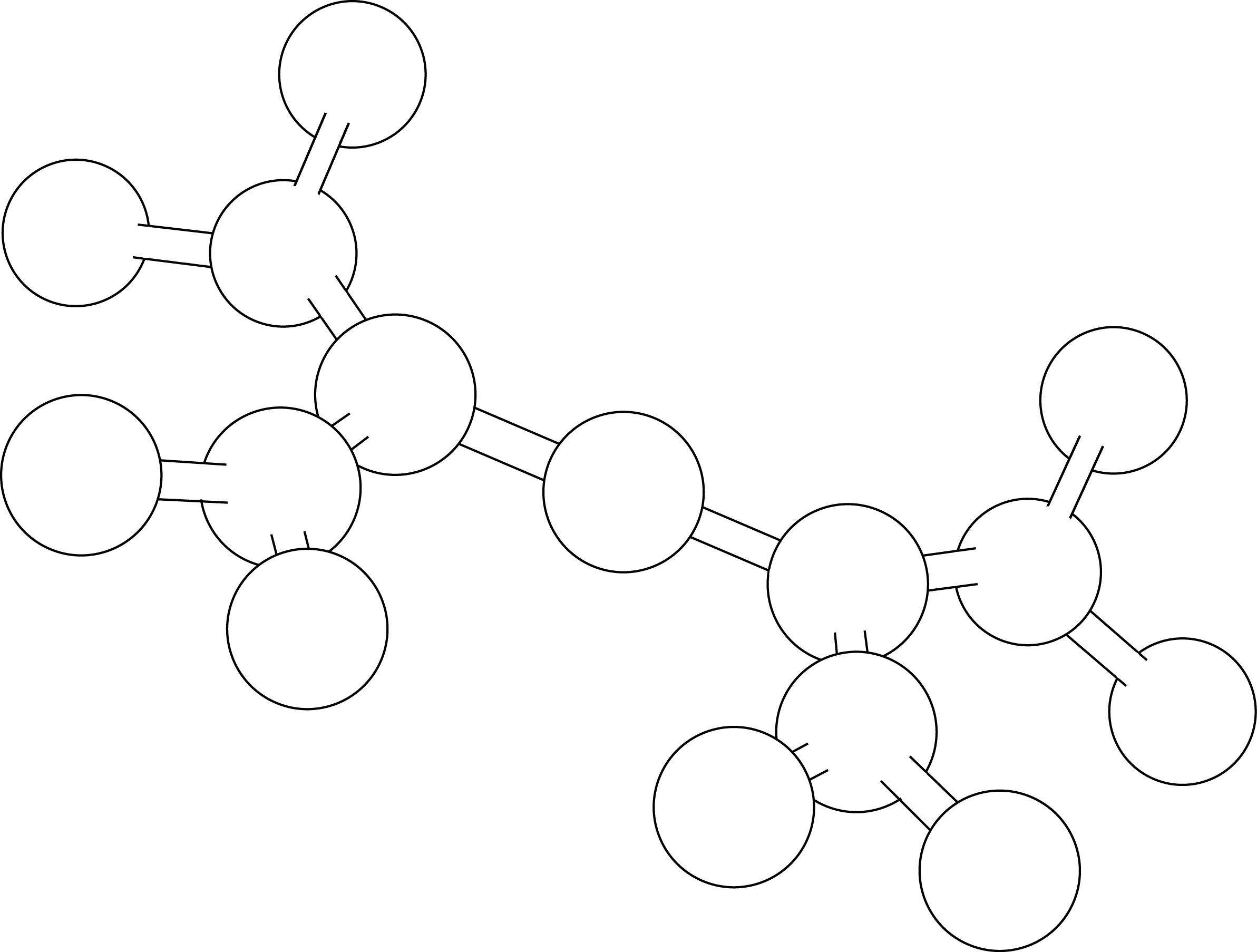}
	\caption{The metric on $S^3$ associated with the full binary tree with $8$ leaves.}\label{figure1}
\end{figure}

In order to prove Theorem A, we construct a special sequence of metrics on the three-sphere using the connected sum procedure of Gromov and Lawson \cite{G-L}. More precisely, given a full binary tree, we associate a spherical region for each node. These regions are subsets of the round three-sphere obtained by removing either one, two or three identical geodesic balls, depending on the vertex degree of the node only. Regions corresponding to neighboring nodes are glued together using a copy of a fixed tube, which is obtained by Gromov-Lawson's method. In Figure \ref{figure1} it is provided a rough depiction of the metric associated with the full binary tree with $8$ leaves.

The lower bounds that we obtain for the widths rely on a combinatorial argument and on the relative isoperimetric inequality. The first key step of our argument is the choice of a special slice $\Sigma_{t_0} = \partial \Omega_{t_0}$, for any fixed sweepout $\{\Sigma_t\}$ of $S^3$ with the metrics that we consider. Then, we induce a $2$-coloring on the nodes of the associated binary tree, where the color of each node depends on the volume of $\Omega_{t_0}$ on the corresponding spherical region. A $2$-coloring of a tree is an assignment of one color, black or white, for each node. Finally, a joint application of a combinatorial tool about $2$-colorings of full binary trees and the relative isoperimetric inequality on some compact three-manifolds with boundary gives us lower bounds on the width.

Liokumovich also used combinatorial arguments to construct Riemannian metrics of large widths on surfaces of small diameter, see \cite{Liokumovich}.

The foundational idea of the Almgren-Pitts min-max theory for the area functional is to achieve the width as the area of a closed minimal surface, possibly disconnected and with multiplicities, see \cite{Pitts} or \cite{C-D} for further details. In our examples of metrics on $S^3$, it is expected that this min-max minimal surface will have multiple components, some of them stable with area strictly less than $4\pi$. These surfaces correspond to the spherical slices of minimum area in the tubes of our construction.     

Brendle, Marques and Neves \cite{BMN} constructed non-spherical metrics with scalar curvature $R\geq 6$ on the hemisphere, which coincide with the standard round metric in a neighborhood of the boundary sphere. These metrics are counterexamples to the Min-Oo conjecture. Our result gives also a setting in which a scalar curvature rigidity result does not hold.

%controlling the width by the sup of the isop. profile
The standard argument to provide a positive lower bound for the width of a Riemannian metric is to consider the supremum of the isoperimetric profile. Given a Riemannian metric $g$ on the three-sphere, its isoperimetric profile is the function $\mathcal I : [0, vol_g(S^3)] \rightarrow \R$ defined by
\begin{equation*}
\mathcal I (v) = \inf \{Area_g(\partial \Omega) : \Omega \subset S^3 \text{ and } vol_g(\Omega)=v\}.
\end{equation*}

In particular, if $\{\Sigma_t\}$ is a sweepout of $(S^3, g)$ and the associated open subsets are $\Omega_t\subset S^3$, then the volumes $vol_g(\Omega_t)$ assume all the values between zero and $vol_g(S^3)$. Then, we have
\begin{equation*}
\sup\{\mathcal I(v) : v \in [0, vol_g(S^3)]\} \leq W(S^3, g).
\end{equation*}
If $g$ has positive Ricci curvature and scalar curvature $R \geq 6$, the $4\pi$ upper bound on the supremum of the left-hand-side of the above expression was previously obtained by Eichmair \cite{Eichmair}. We observed that the supremum of the isoperimetric profiles of the examples that we use also form an unbounded sequence. This claim is also proved via a combinatorial result and provides us a second proof of the content of Theorem A.

\subsection*{Acknowledgments} The results contained in this paper are based partially on the author's Ph.D. thesis under the guidance of Professor Fernando Cod\'a Marques. This work was done while I was visiting him at Princeton University. It is a pleasure to show my gratefulness for his support. 

%%%%%%%%%%%%%%%

\subsection*{Organization} The content of this paper is organized as follows:

In Section \ref{comb-section}, we develop the combinatorial tools that we use to estimate the geometric objects. In Section \ref{GL-section}, we briefly discuss the connected sum procedure due to Gromov and Lawson, and introduce our examples. In Section \ref{width-section}, we prove Theorem A. In Section \ref{isop-section}, we discuss about the isoperimetric profiles of the constructed metrics.

            %%%%%%%%%%
        %%%%%%%%%%%%%%%%%%    
     %%%%%%%%%%%%%%%%%%%%%%%%
   %%%%%%%%%%%%%%%%%%%%%%%%%%%% 
  %%%%%%%%%%%%%%%%%%%%%%%%%%%%%%
 %%%%%%%%%%%%%%%%%%%%%%%%%%%%%%%%%
%%%%%%%%%%%%%%%%%%%%%%%%%%%%%%%%%%%           
%%%%%%%Comeco do trabalho%%%%%%%%%%
%%%%%%%%%%%%%%%%%%%%%%%%%%%%%%%%%%%
 %%%%%%%%%%%%%%%%%%%%%%%%%%%%%%%%%
  %%%%%%%%%%%%%%%%%%%%%%%%%%%%%%
	 %%%%%%%%%%%%%%%%%%%%%%%%%%%% 
	   %%%%%%%%%%%%%%%%%%%%%%%%
		   	%%%%%%%%%%%%%%%%%%  
			     	%%%%%%%%%%

\section{Combinatorial Results}\label{comb-section}

In this section we state and prove our combinatorial results. For each positive integer $m \in \mathbb{N}$, we use $T_m$ to denote the full binary tree, for which all the $2^m$ leaves, nodes of vertex degree one, have depth $m$. Recall that the vertex degree of a node is the number of edges incidents to it. Observe that $T_m$ has $2^{m+1}-1$ nodes, being $2^k$ nodes on the $k$-th level of depth.

We consider $2$-colorings of the nodes of $T_m$. A $2$-coloring of $T_m$ is an assignment of one color, black or white, for each node. Fixed a $2$-coloring, an edge is said to be dichromatic if it connects nodes with different colors. We are interested in $2$-colorings which minimize the quantity of dichromatic edges for a fixed number of black nodes.

\begin{figure}[H]
  \centering
    \caption{A $2$-coloring of $T_2$ with three black nodes and one dichromatic edge.}
    \includegraphics[width=0.3\textwidth]%
    {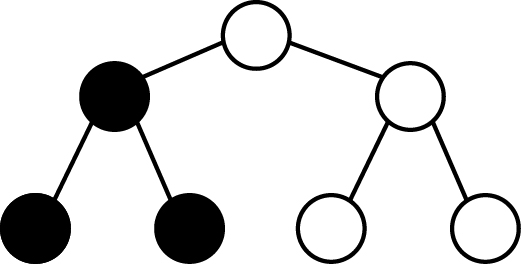}% picture filename
\end{figure}

\subsection{Definition} Let $m, d \in \N$ be positive integers. We define $B_m(d)$ to be the set of values $b$, with $1\leq b \leq 2^{m+1}-1$, for which there exists a $2$-coloring of $T_m$ with $d$ dichromatic edges and $b$ black nodes exactly.

The first statement that we prove in this section is the following upper bound on the size of the sets $B_m(d)$: 

\subsection{Lemma}\label{coloring.nodes}
\textit{For every $m, d \in \N$, it follows that $\# B_m(d) \leq 2^{d} m^{d}$.
}

\begin{proof}
Consider a $2$-coloring $\mathcal C$ of the nodes of $T_m$ which has $d$ dichromatic edges exactly. Let $k$ be the highest level of depth of $T_m$ for which one of those $d$ edges joins a node on the $k$-th level to a node on the $(k+1)$-th level of depth. Choose one of these deeper dichromatic edges and observe that it determines a monochromatic component of $\mathcal C$ which is a copy of $T_{m-(k+1)}$, for some $0 \leq k \leq m-1$. Changing the color of the nodes of this $T_{m-(k+1)}$ yields a $2$-coloring $\mathcal C^{\prime}$ with $d-1$ dichromatic edges and $b$ black nodes, for some $b \in B_m(d-1)$. Since we changed the colors of the nodes in the copy of $T_{m-(k+1)}$ only and they all have the same color on $\mathcal C$, the number of black nodes of $\mathcal C$ is either $b + (2^{m-k} -1)$ or $b - (2^{m-k} - 1)$. Therefore, $\# B_m(d) \leq 2m\cdot \# B_m(d-1)$ and the statement follows inductively and from the fact that $B_m(0) = \{2^{m+1}-1\}$.    
\end{proof}

\subsection{Definition} For $m \in \N$ and $1\leq b \leq 2^{m+1}-1$, let $d_m^{\prime}(b)$ denote the minimum integer $d \in \N$ for which $b \in B_m(d)$. In other words, any $2$-coloring of $T_m$ with $b$ black nodes exactly has at least $d_m^{\prime}(b)$ dichromatic edges.

As a consequence of the above estimate we prove a qualitative result that guarantees that there is no uniform bound on the values $d_m^{\prime}(b)$.

\subsection{Proposition}\label{coloring.nodes.2}
\textit{There exist integers $b(m)$, with $1\leq b(m) \leq 2^{m+1}-1$ and such that $\{d_m^{\prime}(b(m))\}_{m\in \N}$ is an unbounded sequence.
}

\begin{proof}
Suppose, by contradiction, there exists $D \in \N$ such that $d_m^{\prime}(b) \leq D$, for all $m\in \N$ and $1\leq b \leq 2^{m+1}-1$. Then, each such $b$ belongs to $B_m(d)$, for some $d \leq D$. By Lemma \ref{coloring.nodes}, we have
\begin{equation*}
\# \bigg( \bigcup_{d \leq D} B_m(d) \bigg) \leq \sum_{d \leq D} 2^{d} m^{d}.
\end{equation*}
But, for fixed $m\in \N$, there are $2^{m+1}-1$ possible values for $b$, which can not be controlled by the polynomial right hand side of the above expression. 
\end{proof}

Since the vertex degrees of the nodes of the considered trees do not exceed $3$, we can use Proposition \ref{coloring.nodes.2} to prove the following:

\subsection{Corollary}\label{disjoint-dich-edges}
\textit{Given $k \in \N$, there exist $m \in \N$ and $b(m) < 2^{m+1}-1$, such that any $2$-coloring of $T_m$ with $b(m)$ black nodes exactly has at least $k$ pairs of neighboring nodes with different colors. Moreover, for $1\leq b \leq 2^{m+1}-1$, any $2$-coloring of $T_m$ with $b$ black nodes exactly has at least $(k-|b-b(m)|)/5$ pairwise disjoint pairs of neighboring nodes with different colors.
}

\begin{proof}
Indeed, let $m$ and $b(m)$ be such that $d_m^{\prime}(b(m))\geq k$. This choice is allowed by Proposition \ref{coloring.nodes.2}. Then, any $2$-coloring of $T_m$ with $b(m)$ black nodes exactly has at least $k$ pairs of neighboring nodes with different colors.

For $1\leq b \leq 2^{m+1}-1$, we can estimate $d_m^{\prime}(b)$ using the formula:
\begin{equation}\label{eq-A}
|d_m^{\prime}(t) - d_m^{\prime}(s)| \leq |t-s|.
\end{equation}
To prove this relation, we need to verify $|d_m^{\prime}(t) - d_m^{\prime}(t+1)| \leq 1$ only. Consider a $2$-coloring of $T_m$ with $t+1$ black nodes and $d_m^{\prime}(t+1)$ dichromatic edges exactly. Let $N$ be a black node of this coloring with the property that no other black node lives in a level deeper than its level of depth. Changing the color of $N$ to white we obtain a $2$-coloring with exactly $t$ black nodes and at most $d_m^{\prime}(t+1)+1$ dichromatic edges. This implies that $d_m^{\prime}(t)\leq d_m^{\prime}(t+1)+1$. Analogously, we obtain $d_m^{\prime}(t+1)\leq d_m^{\prime}(t)+1$, and we are done with the proof of expression (\ref{eq-A}).

The choice of $m$ and $b(m)$, together with equation (\ref{eq-A}) gives us that $d_m^{\prime}(b) \geq k-|b-b(m)|$. Observe that each pair of neighboring nodes has a common node with four other pairs of neighboring nodes at most. This allows us to conclude that any $2$-coloring of $T_m$ with $b$ black nodes exactly has at least $(k-|b-b(m)|)/5$ pairwise disjoint pairs of neighboring nodes with different colors. And this concludes the proof of the corollary.
\end{proof}

%%%%%%%%%%%%%%

Corollary \ref{disjoint-dich-edges} is key in the proof of the lower bound that we provide for the supremum of the isoperimetric profiles of the Riemannian metrics that we construct in the next section. This is done in Section \ref{isop-section}.

The rest of this section is devoted to the discussion of an interesting quantitative statement related to the previous results. It also can be applied to provide estimates for the widths of our examples.

%%%%%%%%%%%%dichromatic value

\subsection{Definition} We define the \textit{dichromatic value of $t$ leaves in $T_m$} as the least number of dichromatic edges of a $2$-coloring of $T_m$ with $t$ black leaves exactly. We denote this number by $d_m(t)$.  

The following statement about dichromatic values is the analogous of formula (\ref{eq-A}) for a fixed number of black leaves. We omit its proof here.

\subsection{Lemma}\label{di-value}
\textit{ $|d_m(t) - d_m(s)| \leq |s-t|$.
}

%%%%%%%%%%%%Combinatorial theorem

\subsection{Theorem}\label{coloring.leaves}
\textit{ For each integer $m > 1$, let
\[ a(m) = \left\{ 
  \begin{array}{l l}
   1 + 2 + 2^3 + \ldots + 2^{m-2}, & \quad \text{if } m \text{ is odd}\\
   1 + 2^2 + 2^4 + \ldots + 2^{m-2}, & \quad \text{if } m \text{ is even}.
  \end{array} \right.\]
Then, $d_m(a(m)) \geq \left\lceil \frac{m}{2} \right\rceil$.
}

%%%%%%%%%%%%%Proof of combinatorial theorem

\begin{proof}
The proof is by induction. Define $a(1) = 1$. The initial cases, $m=1$ and $2$, are very simple. Suppose that the statement is true for $m-1$ and $m-2$. Let $\C$ be a $2$-coloring of $T_m$ with $a=a(m)$ black leaves exactly. 

Let us count the number of different types of nodes in the $(m-1)$-th level of depth of $T_m$. Use $\alpha$ to denote the number of nodes which have two neighbor leaves of different colors on $\C$ and, similarly, let $\beta$ be the number of nodes which have two black neighbor leaves. Then, we can write the number of black leaves as $a(m) = \alpha +2\cdot \beta$. In particular, this expression implies that $\alpha$ is an odd number and

\begin{equation}\label{eq1}
\beta = \frac{1-\alpha}{2} + \frac{a(m)-1}{2} = \frac{1-\alpha}{2} + a(m-1) - m^{\prime},
\end{equation}
where $m^{\prime} = 1$ if $m$ is even and $m^{\prime} = 0$ if $m$ is odd. Indeed, it is easily seen that $a(m)-1=2\cdot (a(m-1) - m^{\prime})$, for every $m\geq 3$.

Next, we induce a $2$-coloring $\C ^{\prime}$ on the nodes of $T_{m-1}$. Let $\T_{m-2} \subset \T_{m-1} \subset T_m$ be such that $T_m$ minus $\T_{m-1}$ is the set of leaves of $T_m$, and $\T_{m-1}$ minus $\T_{m-2}$ is the set of leaves of $\T_{m-1}$. We begin to define $\C^{\prime}$ on $T_{m-1}$, identified with $\T_{m-1}$, asking $\C^{\prime}$ to be equal to $\C$ on $\T_{m-2}$ and on the $\alpha$ leaves of $\T_{m-1}$ which have neighbor leaves of different colors on $\C$. The $\beta$ leaves of $\T_{m-1}$ which have two black neighbor leaves of $\C$ on $T_m$ are colored black on $\C^{\prime}$. The remaining leaves of $\T_{m-1}$ have two white neighbor leaves of $\C$ on $T_m$ and receive the color white on $\C^{\prime}$. The important properties of $\C^{\prime}$ are:

\begin{enumerate}
\item[(i)] The number of black leaves of $\C^{\prime}$ on $\T_{m-1}$ is greater than or equal to $\beta$ and at most $\alpha + \beta$; 

\item[(ii)] The number of dichromatic edges of $\C$ is at least $\alpha$ plus the number of dichromatic edges of $\C^{\prime}$. 
\end{enumerate}

The first of these properties follows directly from the construction. To prove the second, we begin by observing that the dichromatic edges of $\C^{\prime}$ that are in $\T_{m-2}$ are, automatically, dichromatic edges of $\C$ on $T_m$. Then, we analyze cases to deal with the dichromatic edges of $\C^{\prime}$ that use leaves of $\T_{m-1}$. We omit this simple analysis.   

Let us use $t$ to denote the number of black leaves of $\C^{\prime}$ on $\T_{m-1}$. From (\ref{eq1}) and property (i) above, we conclude that
\begin{equation}\label{eq2}
\frac{1-\alpha}{2} - m^{\prime} \leq t - a(m-1) \leq  \frac{1+\alpha}{2} - m^{\prime}.
\end{equation}
By Lemma \ref{di-value} and the induction hypothesis, we have
\begin{equation}\label{eq-3}
d_{m-1}(t) \geq d_{m-1}(a(m-1)) - \frac{1+\alpha}{2} \geq \left\lceil \frac{m-1}{2} \right\rceil - \frac{1+\alpha}{2}.
\end{equation}
By definition, $\C^{\prime}$ has at least $d_{m-1}(t)$ dichromatic edges. This implies, together with property (ii) and equation (\ref{eq-3}), that
\begin{equation}\label{eq-4}
\#\{\text{dichromatic edges of } \C \} \geq \left\lceil \frac{m-1}{2} \right\rceil + \frac{\alpha-1}{2}.
\end{equation}  

Since $\alpha$ is odd, if $m$ is even we have that the number of dichromatic edges of $\C$ is at least $\left\lceil (m-1)/2 \right\rceil = \left\lceil m/2 \right\rceil$ and we are done. Otherwise, $m$ is odd and equation (\ref{eq-4}) provides us
\begin{equation*}
\#\{\text{dichromatic edges of } \C \} \geq \left\lceil \frac{m}{2} \right\rceil - 1 + \frac{\alpha-1}{2} = \left\lceil \frac{m}{2} \right\rceil + \frac{\alpha-3}{2}.
\end{equation*}
If $\alpha \geq 3$, the induction process ends and we are done. 

From now on, we suppose that $m$ is odd and $\alpha = 1$. The arguments are going to be similar to the previous one, but one level above on $T_m$. Recall that $\alpha = 1$ means that there is a unique node of $T_m$ which is neighbor of leaves with different colors. Observe that each node on the $(m-2)$-th level of $T_m$, leaf of $\T_{m-2}$, is associated with four leaves of $T_m$, the leaves at edge-distance two. A unique node on this level is associated with an odd number ($1$ or $3$) of black leaves of $\C$, because $\alpha = 1$. We denote this node and odd number by $N$ and $\varphi$, respectively. 

Let us use $\theta$ and $\gamma$ to denote the number of nodes on the $(m-2)$-th level of $T_m$ which are associated, respectively, with two and four black leaves of $\C$. Then, we can write the number of black leaves of $\C$ as $a(m) = \varphi + 2 \theta + 4 \gamma$.

In particular, this expression implies that $\varphi + 2\theta \equiv 3$ mod $4$ and

\begin{equation*}
\gamma = - \frac{1 + \varphi + 2 \theta}{4} + a(m-2).
\end{equation*}

Next, we induce a $2$-coloring $\C ^{\pprime}$ on the nodes of $T_{m-2}$. We follow the same $\T_{m-2} \subset T_m$ notation that we introduced before. And now, we still need the analogous $\T_{m-3} \subset \T_{m-2}$. We choose $\C^{\pprime}$ to coincide with $\C$ on the nodes of $\T_{m-3}$ and on the $\theta$ leaves of $\T_{m-2}$ which are associated with two black leaves of $\C$. The $\gamma$ leaves of $\T_{m-2}$ which are associated with four black leaves of $\C$ are colored black. The node $N$ is colored white if $\varphi = 1$ and black if $\varphi = 3$. The leaves of $\T_{m-2}$ which remain uncolored receive the color white. As in the previous step, the important properties of $\C^{\pprime}$ are:

\begin{enumerate}
\item[(I)] The number of black leaves of $\C^{\pprime}$ is between the values $(\varphi -1)/2 + \gamma$ and $(\varphi -1)/2 + \gamma + \theta$; 

\item[(II)] The number of dichromatic edges of $\C$ is at least $1+ \theta$ plus the number of dichromatic edges of $\C^{\pprime}$. 
\end{enumerate}

By the same reasoning that we used to obtain equation (\ref{eq-4}), we have

\begin{equation*}
\#\{\text{dichromatic edges of } \C \} \geq \left\lceil \frac{m-2}{2} \right\rceil + \frac{1+ \varphi + 2\theta}{4}.
\end{equation*}
Using that $m$ is odd, $\varphi = 1$ or $3$ and $\varphi + 2\theta \equiv 3$ mod $4$, we conclude that the right hand side of the above expression is greater than or equal to $\left\lceil m/2 \right\rceil$. This finishes the induction step.
\end{proof}

%%%%%%%%%%%%%%%
An easy consequence of Theorem \ref{coloring.leaves} is the following:  

\subsection{Corollary}\label{disjoint-dich-edges2}
\textit{Given $m \in \N$ there exists $a(m) \in \N$, $1\leq a(m) \leq 2^m$, such that any $2$-coloring of $T_m$ with $a(m)$ black leaves exactly has at least $(\left\lceil m/2\right\rceil)/5$ pairwise disjoint pairs of neighboring nodes with different colors.
}

%%%%%%%%%%%%%%%%%%%%%%%

\section{Constructing the examples}\label{GL-section}

In this section, we introduce our examples. We begin with a brief discussion about the Gromov-Lawson metrics of positive scalar curvature. Then, for each full binary tree we construct an associated metric of scalar curvature greater than or equal to $6$ on the three-sphere.

\subsection{Gromov-Lawson metrics} 

Gromov and Lawson developed a method that is adequate to perform connected sums of manifolds with positive scalar curvature, see \cite{G-L}. They proved the following statement:  

Let $(M^n, g)$ be a Riemannian manifold of positive scalar curvature. Given $p \in M$, $\{e_1,\ldots, e_n\} \subset T_p M$ an orthonormal basis, and
$r_0>0$, it is possible to define a positive scalar curvature metric $g^{\prime}$ on the punctured geodesic ball $B(p,r_0)-\{p\}$ that coincides with $g$ near the boundary $\partial B(p,r_0)$, and such that $(B(p,r_1)-\{p\}, g^{\prime})$ is isometric to a half-cylinder for some $r_1>0$.

\subsubsection{Remark:}\label{size of tube} Since $(B(p,r_1)-\{p\}, g^{\prime})$ is isometric to a half-cylinder, there exists $0 < r <r_1$ so that the $g^{\prime}$ volume of $B(p,r_0) - B(p,r)$ is bigger than one half of the $g$ volume of the removed geodesic ball $B(p,r_0)$.

\subsection{Fundamental blocks}\label{blocks}

In order to construct our examples, we use the above metrics to build three types of fundamental blocks. 

The first type is obtained from the above construction using the standard round metric on $M = S^3$, any $p \in S^3$ and orthonormal basis and $r_0 = 1$. We use $S_1 = (S^3 - B(p, r), g^{\prime})$ to denote this block, where $0< r < r_1 < r_0 = 1$ is chosen as in remark \ref{size of tube}. Observe that $S_1$ is a manifold with boundary, has positive scalar curvature and it has a product metric near the boundary two-sphere $\partial S_1$. To obtain the second fundamental block, $S_2$, we perform the same steps with $M = S_1$ and choosing the new removed geodesic ball to be antipodally symmetric to $B(p,1)$. Finally, the third block, $S_3$, is obtained from $S^3$ after three application of Gromov-Lawson procedure by removing three disjoint geodesic balls $B(p_i,1)$, $i= 1, 2$ and $3$.

Summarizing, the fundamental blocks $S_1, S_2$ and $S_3$ are obtained from the standard three-sphere by removing one, two or three geodesic balls, respectively, and attaching a copy of a fixed piece for each removed ball. Each fundamental block has a product metric near its boundary spheres. Up to a re-scaling, admit that they have scalar curvature greater than or equal to $6$. Also, by remark \ref{size of tube}, we can suppose that the attached piece has volume bigger than one half of the volume of the removed balls.

\subsection{A metric on $S^3$ associated with $T_m$}

For each full binary tree $T_m$, we use the fundamental blocks to construct an associated metric on $S^3$.

In our examples, each node of $T_m$ will be associated to a fundamental block. Following the notation of subsection \ref{blocks}, for a node of degree $k$ we associate a fundamental block of type $S_k$. We connect the blocks which correspond to neighboring nodes of $T_m$ by identifying one boundary sphere of the first to a boundary sphere of the second with reverse orientations. After performing all identifications, we obtain a metric on $S^3$, which is denoted by $g_m$ and has scalar curvature $R\geq 6$. 

The metric $g_m$ decomposes $S^3$ in $2^{m+1}-1$ disjoint closed regions, which are isometric to the standard three-sphere with either one, two or three identical disjoint geodesic balls removed, and $2^{m+1}-2$ connecting tubes. Moreover, by construction, the tubes are isometric to each other and their volume is greater than the volume of one of the removed geodesic balls.

%%%%%%%%%%%%%%%%%%%%%%%
%%%%%%%%%%%%%%%%%%%%%%%

\section{Lower bounds on the width}\label{width-section}

In the introduction, we briefly defined the notions of sweepouts and width of Riemannian metrics on $S^3$. We begin this section by recalling what these interesting geometrical objects are. Then, we prove that the widths of the metrics $g_m$ constructed in Section \ref{GL-section} converge to infinity.

We use $I = [0,1]$ to denote the closed unit interval on the real line. The $2$-dimensional Hausdorff measure on $S^3$ induced by a Riemannian metric $g$ is denoted by $\H^2$. Let us remember the definition of sweepouts.

A \textit{sweepout} of $(S^3,g)$ is a family $\{\Sigma_t\}_{t\in I}$ of smooth $2$-spheres, which are boundaries of open sets $\Sigma_t = \partial \Omega_t$ such that:
\begin{enumerate}
\item $\Sigma_t$ varies smoothly in $(0,1)$;
\item $\Omega_0 = \varnothing$ and $\Omega_1 = S^3$;
\item $\Sigma_t$ converges to $\Sigma_{\tau}$, in the Haursdorff topology, as $t\rightarrow \tau$;
\item $\H^2(\Sigma_t)$ is a continuous function of $t \in I$.
\end{enumerate}

Let $\Lambda$ be a set of sweepouts of $(S^3, g)$. It is said to be saturated if given a map $\phi \in C^{\infty}(I\times S^3, S^3)$ such that $\phi(t, \cdot)$ are diffeomorphism of $S^3$, all of which isotopic to the identity, and a sweepout $\{\Sigma_t\}_{t\in I} \in \Lambda$, we have $\{\phi(t,\Sigma_t)\}_{t\in I} \in \Lambda$. The \textit{width of $(S^3,g)$ associated with $\Lambda$} is the following min-max invariant:
\begin{equation*} 
W(S^3, g, \Lambda) = \inf_{\{\Sigma_t\} \in \Lambda} \ \ \max_{t\in [0,1]} \H^2(\Sigma_t).
\end{equation*}

From now on, we fix a saturated set of sweepouts $\Lambda$. The main result of this work is the following:

\subsection{Theorem}\label{thm1}
\textit{The sequence $\{g_m\}_{m\in \N}$ of Riemannian metrics on $S^3$ that we constructed satisfies:
\begin{equation*}
\lim_{m\rightarrow \infty} W(S^3, g_m, \Lambda) = + \infty.
\end{equation*} 
}

\begin{proof}
Recall that $g_m$ is a metric on $S^3$ which is related to the full binary tree $T_m$. Also, there are $2^{m+1}-1$ disjoint closed subsets of $(S^3, g_m)$ isometric to the standard round metric on the three-sphere with either one, two or three identical disjoint balls removed. These spherical regions are associated to nodes of $T_m$ and two of them are glued to each other if, and only if, their corresponding nodes are neighbors in $T_m$. It is also important to recall that there is a tube connecting such neighboring regions, all of which isometric to each other. For convenience, let us denote these spherical regions by:
\begin{itemize}
\item $L_i$ if it corresponds to a leaf of $T_m$;
\item $B$ if it corresponds to the node of degree $2$;
\item $A_j$ if it corresponds to a node of degree $3$. 
\end{itemize}

After connecting any two neighboring spherical regions, we obtain a region $\mathcal A$ which is isometric to one of the domains depicted in Figure \ref{figure2}.

\begin{figure}[H]
  \centering
    \includegraphics[width=1.0\textwidth]{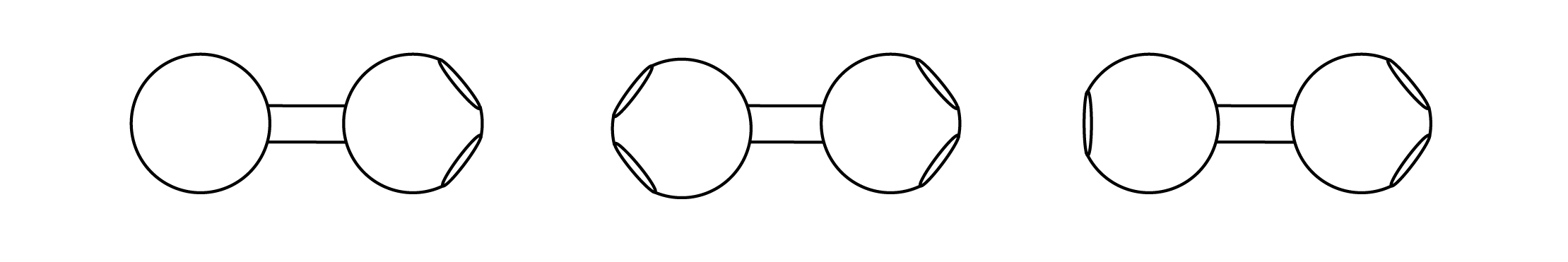}
\caption{The possible regions that we obtain after gluing together two neighboring spherical regions.}\label{figure2}
\end{figure}

In the figure, the first region was obtained by gluing a spherical region of type $L_i$ to its only $A_j$ neighbor. The others are obtained by gluing either two $A_j$ regions or the $B$ region to one its two $A_j$ neighbors. 

Choose $\alpha> 0$ such that $2\alpha$ is strictly less then the volume of one $A_j$. Then, $0< \alpha < vol(\mathcal A) -\alpha < vol(\mathcal A)$, for any of the possible $\mathcal A$'s. By the relative isoperimetric inequality, there exists $C>0$ such that for any open subset $\Omega \subset \mathcal A$ of finite perimeter and $\alpha \leq vol(\Omega) \leq vol(\mathcal A) -\alpha$, we have
\begin{equation}\label{isop-inequality}
\H^2(\partial \Omega \cap int(\mathcal A)) \geq C.
\end{equation}  
Moreover, since we have three possible isometric types of $\mathcal A$'s only, we can suppose that this constant does not depend on the type of $\mathcal A$.

Let $\{\Sigma_t\}_{t\in I}$ be a sweepout of $(S^3, g_m)$. Consider the associated open sweepout $\{\Omega_t\}$, for which $\Sigma_t = \partial \Omega_t$. Let $a(m)$ be the integer provided by Corollary \ref{disjoint-dich-edges2}. Choose the least $t_0 \in I$ for which we have $vol(\Omega_{t_0} \cap L_i) \geq \alpha$, for at least $a(m)$ values of $i\in \{1, 2, 3, \ldots, 2^m\}$. 

Observe that, at most $a(m)-1$ of the $L_i$'s can satisfy $vol(\Omega_{t_0}\cap L_i) > \alpha$. Up to a reordering of their indices, suppose that $vol(\Omega_{t_0} \cap L_i) \geq \alpha$, for $i = 1, 2, \ldots, a(m)$, and $vol(\Omega_{t_0} \cap L_i) \leq \alpha$, otherwise. 

On $T_m$, consider the $2$-coloring defined in the following way: the leaves associated to $L_1, \ldots, L_{a(m)}$ are colored black, the other leaves are colored white and the nodes which are not leaves are colored black if, and only if, the volume of $\Omega_{t_0}$ inside the corresponding spherical region is greater than or equal to $\alpha$. This $2$-coloring of $T_m$ has exactly $a(m)$ black leaves.

By Corollary \ref{disjoint-dich-edges2}, the constructed coloring has at least $m/10$ pairwise disjoint pairs of neighboring nodes with different colors. Observe that each such pair gives one $\mathcal A$ type region for which we have
\begin{equation}\label{volume of omega}
\alpha \leq vol(\Omega_{t_0}\cap \mathcal A) \leq vol(\mathcal A) -\alpha. 
\end{equation}
This follows because we chose $\alpha$ in such a way that the volume of our spherical regions are greater than $2\alpha$. By equation (\ref{isop-inequality}) and (\ref{volume of omega}) we conclude that $\H^2(\partial \Omega_{t_0}\cap int(\mathcal A)) \geq C$. Since this holds for $m/10$ pairwise disjoint $\mathcal A$ type regions, we have $\H^2(\Sigma_{t_0}) \geq C\cdot m/10$. This concludes our argument.
\end{proof}

%%%%%%%%%%%%%%%%%%%%%%%%%%%%
%%%%%%%%%%%%%%%%%%%%%%%%%%%%

\section{On the isoperimetric profiles of the metrics $g_m$} \label{isop-section}

In this section, we discuss the fact that isoperimetric profiles $\mathcal I_{m}$ of the metrics $g_m$ are not uniformly bounded. This part also relies on a combinatorial argument, we use Corollary \ref{disjoint-dich-edges}. The idea is to decompose $S^3$ into $2^{m+1}-1$ pieces of identical $g_m$ volumes, all of which being the union of one spherical region ($L_i$, $B$ or $A_j$) with a portion of their neighboring tubes. 

\begin{figure}[H]
 \centering
  \includegraphics[width=.9\textwidth]{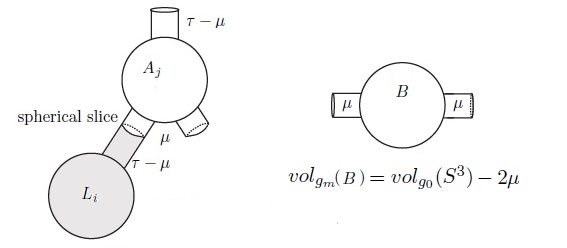}
\caption{The balanced regions.}\label{figure4}
\end{figure}

%About the volume of the tube
This decomposition of $S^3$ by balanced pieces is only possible because we chose the tube large enough to have $g_m$ volume greater than the spherical volume $\mu$ of the removed geodesic balls. This allows us to decompose $S^3$ into $2^{m+1}-2$ pieces with volume $vol_{g_0}(S^3) + \tau - 2\mu$ and one piece with volume $vol_{g_0}(S^3)$, and with the other desired properties, where $g_0$ is the standard round metric on $S^3$ and $\tau$ is the volume of the gluing tube. The piece with volume $vol_{g_0}(S^3)$ is the one related to the only node of degree $2$ in $T_m$. For each pair of neighboring spherical regions, the boundary of the associated balanced regions has exactly one component in the connecting tube. This component is a spherical slice which splits the volume $\tau$ of the tube as $\tau - \mu$ plus $\mu$, as depicted in Figure \ref{figure4}. 

There are three types of balanced regions, depending on the number of boundary components, all of which we denote by $\mathcal M$. In this part, we use the isoperimetric inequality in its full generality: there exists $C>0$ so that
\begin{equation}\label{full.isop.ineq}
\min \{\H^3(\Omega), \H^3(\mathcal M - \Omega)\}^{2/3} \leq C\cdot \H^2(\partial \Omega \cap int(\mathcal M)),
\end{equation}
for every $\Omega \subset \mathcal M$. Since we have three types of $\mathcal M$ regions only, we suppose that $C>0$ associated to $\mathcal M$ does not depend on its type.

Suppose, by contradiction, that there exists $L>0$ such that $\mathcal I_m(v) \leq L$, for every $m \in \N$ and $v \in [0,vol_{g_m}(S^3)]$.

For any $k \in \N$, let $m, b(m) \in \N$ be the integers provided by Corollary \ref{disjoint-dich-edges}. Take $v(m) = b(m)\cdot (vol_{g_0}(S^3) + \tau - 2\mu)$ and let $\Omega \subset S^3$ be a subset with $vol_{g_m}(\Omega) = v(m)$ and such that $\H^2(\partial \Omega) \leq L$. Observe that
\begin{equation*}
\sum_{\mathcal M} \H^2(\partial \Omega \cap int(\mathcal M)) \leq \H^2(\partial \Omega) \leq L.
\end{equation*}
Using the relative isoperimetric inequality, equation (\ref{full.isop.ineq}), we obtain:
\begin{equation*}
\sum_{\mathcal M} \min \{\H^3(\Omega\cap \mathcal M), \H^3(\mathcal M - \Omega)\}^{2/3} \leq C\cdot L.
\end{equation*}
Which implies that
\begin{equation}\label{eq-A1}
\sum_{\mathcal M} \min \{\H^3(\Omega\cap \mathcal M), \H^3(\mathcal M - \Omega)\}\leq C_1,
\end{equation}
where $C_1 = (C\cdot L)^{3/2}$. Let $\mathcal M_1$ be the set of the $\mathcal M$ type regions for which $2\cdot \H^3(\Omega \cap \mathcal M) < \H^3 (\mathcal M)$. Similarly, the $\mathcal M$ regions satisfying the opposite inequality compose $\mathcal M_2$. Equation (\ref{eq-A1}) implies
\begin{equation}
\sum_{\mathcal M \in \mathcal M_1} \H^3(\Omega\cap \mathcal M) + \sum_{\mathcal M \in \mathcal M_2} \H^3(\mathcal M - \Omega) \leq C_1.  
\end{equation}
Using that $\H^3(\mathcal M - \Omega) = \H^3(\mathcal M) - \H^3(\Omega \cap \mathcal M)$ and the fact that
\begin{equation}
\sum_{\mathcal M \in \mathcal M_1} \H^3(\Omega\cap \mathcal M) + \sum_{\mathcal M \in \mathcal M_2} \H^3(\Omega\cap \mathcal M)  = vol_{g_m}(\Omega) = v(m), 
\end{equation}
we easily conclude
\begin{equation}\label{eq-A2}
\bigg| v(m) - \sum_{\mathcal M \in \mathcal M_2} \H^3 (\mathcal M)\bigg| \leq C_1.
\end{equation}
By the choice of $v(m)$, equation (\ref{eq-A2}) implies that $|\# \mathcal M_2 - b(m)| \leq C_2$, where $C_2 = (C_1 + |\tau - 2 \mu|)/(vol_{g_0}(S^3) + \tau - 2\mu)$ is a uniform constant.

Consider the $2$-coloring $\mathcal C$ of $T_m$ whose black nodes are those associated with the balanced regions in $\mathcal M_2$. Then, $\mathcal C$ has $\# \mathcal M_2$ black nodes exactly. By Corollary \ref{disjoint-dich-edges}, there are at least $(k - |\# \mathcal M_2 - b(m)|)/5$ pairwise disjoint pairs of neighboring nodes with different colors. By a reasoning similar to the one that we used in the end of Section \ref{width-section}, we have that
\begin{equation}
\H^2(\partial \Omega) \geq C_3 \cdot (k - |\# \mathcal M_2 - b(m)|)/5,
\end{equation} 
for some $C_3>0$, which does not depend on $m$. Recalling that $\H^2(\partial \Omega) \leq L$ and $|\# \mathcal M_2 - b(m)| \leq C_2$, the above expression provides a uniform upper bound on $k$, which is arbitrary. This is a contradiction and we are done.

\end{document}